\documentclass[12pt,draft]{amsart}
\usepackage{a4wide}
\usepackage{amssymb}

\newtheorem{theorem}{Theorem}[section]
\newtheorem{thm}[theorem]{Theorem}
\newtheorem{lemma}[theorem]{Lemma}

\theoremstyle{definition}

\newtheorem{example}[theorem]{Example}

\theoremstyle{remark}

\numberwithin{equation}{section}

\newcommand{\DD}{{\mathbb D}}

\newcommand{\RR}{{\mathbb R}}
\newcommand{\CC}{{\mathbb C}}

\newcommand{\TT}{{\mathbb T}}
\newcommand{\eps}{\varepsilon}

\renewcommand{\phi}{\varphi}

\begin{document}

\title{Shcherbina's theorem for finely holomorphic functions}

\author{Armen Edigarian}

\address{Institute of Mathematics, Jagiellonian University,
\L ojasiewicza 6/2117, 30-348 Krak\'ow, Poland}
\email{armen.edigarian@uj.edu.pl}
\thanks{ }
\author{Jan Wiegerinck}

\address{KdV Institute for Mathematics, University of Amsterdam, Plantage Muidergracht,
24, 1018 TV, Amsterdam, The Netherlands}

\email{j.j.o.o.wiegerinck@uva.nl}






\begin{abstract} We prove an analogue of Sadullaev's theorem concerning the size of the set where a maximal totally real manifold $M$ can meet a pluripolar set. $M$ has to be of class $C^1$ only. This readily leads to a version of Shcherbina's theorem for $C^1$ functions $f$ that are defined in a neighborhood of certain compact sets $K\subset\CC$. If the graph $\Gamma_f(K)$ is
pluripolar, then $\frac{\partial f}{\partial\bar z}=0$ in the closure of the fine interior of $K$.
\end{abstract}

\maketitle


\section{Introduction}
Let $P\subset\CC^n$ be any subset. We say that $P$ is {\em pluripolar} if there exists a plurisubharmonic function $u$ on $\CC^n$, $u\not\equiv-\infty$, such that $P\subset\{u=-\infty\}$.
It is well known that this global definition is equivalent to the local definition (see \cite{Jos}). 

A.~Sadullaev (see \cite{Sadullaev}) proved the following result
\begin{theorem}\label{thm:pl} Let $P\subset\CC^n$ be a pluripolar set and let $M$ be a maximal totally real submanifold 
of class $C^3$ in some domain $D\subset\CC^n$. Then $P\cap M$ has zero measure on $M$.
\end{theorem}
For the definition of a maximal totally real submanifold see Section 3.
B.~Coupet (see \cite{Coupet}) proved Theorem~\ref{thm:pl} for $M$ of class $C^2$. We are interested in a similar
problem but for manifolds $M$ of class $C^1$. Before we present our result, we need the following definition.
Let $L\subset\CC^n$ be any subset and let $z_0\in\CC^n$ be a point. We say that $L$ is thin at $z_0$
if there exists a neighborhood $U$ of $z_0$ in $\CC^n$ and a negative plurisubharmonic function $u$ on $U$
such that $u\le-1$ on $L\setminus\{z_0\}\cap U$ and $u(z_0)>-\frac14$.


\begin{theorem}\label{thm:main1} Let $P\subset\CC^n$ be a pluripolar set and let $M\subset\CC^n$ be a $C^1$ maximal 
totally real manifold. Then $M\setminus P$ is not thin at any point of $M$.
\end{theorem}

The primary  motivation of our paper was the following question: suppose that $E\subset\CC$ is any subset and
$f:E\to\CC$ is any function. What is the relation between pluripotential properties of 
the graph of $f$ over $E$, i.e.,
\begin{equation}
\Gamma_f(E)=\{(z,f(z)): z\in E\},
\end{equation}
and analytic properties of $f$? To be more precise, recall that a set $P\subset\CC^n$ is called {\em pluripolar}
if there exists a plurisubharmonic function $u$ on $\CC^n$, $u\not\equiv-\infty$, such that $P\subset\{u=-\infty\}$.
It is well known that this global definition is equivalent to the local definition (see \cite{Jos}). In particular,
if $E$ is an open set and $f$ is holomorphic on $E$, then $\Gamma_f(E)$ is pluripolar in $\CC^2$. If $E$ is not an open
set (for example, if $E$ is a compact set without interior points) then the situation is very complicated, see e.g.~\cite{Coman-Levenberg-Poletsky} for $E$ equal to the unit circle and $f$ a quasianalytic function, or
\cite{Edi-Marz-Wie} for $E$ a fine domain and $f$ a finely holomorphic function, cf.~Section \ref{fine} below.

The inverse problem, i.e.,  to deduce some analytic properties
of $f$ from the pluripolarity of $\Gamma_f(E)$, seems to be even more difficult.
 Recently, N. Shcherbina, \cite{Shch}, proved the following result, which was conjectured by T.~Nishino.
\begin{theorem}\label{thm:1} Let $D\subset\CC$ be a domain and let $f:D\to\CC$ be a continuous function.
Assume that $\Gamma_f(D)$ is pluripolar in
$\CC^2$. Then $f$ is holomorphic on $D$.
\end{theorem}
Shcherbina also mentioned (see Remark on page 204 in \cite{Shch}) that one can
prove Theorem~\ref{thm:1} for a $C^1$-function $f$ using Bishop's technique (see \cite{Bishop}, c.f. also \cite{K-C}). The assumption that $f$ is not holomorphic would imply that $M=\Gamma_f(D)$ is a totally real manifold, to which a family of analytic discs can be attached, which eventually leads to a contradiction.

What if one drops the assumption that $f$ is defined on a domain? Using results of Coupet (see \cite{Coupet}), which are based on the approach of S.~Pinchuk (see
\cite{Pinchuk}) and A.~Sadullaev (see \cite{Sadullaev}), T.~Edlund proved in his thesis \cite{Edlund} the following result for compact sets. Actually, Edlund states less, but Theorem \ref{thm:2} follows immediately from his proof.
\begin{theorem}\label{thm:2}
Let $K\subset\CC$ be a compact set and let $f:\CC\to\CC$ be a $C^2$-function such that 
the graph of $f$ over $K$ is pluripolar in $\CC^2$.
Put 
\begin{equation*}
S=\{z\in K:\text{ for any }\epsilon>0\text{ the set }K\cap \DD(z,\epsilon)
\text{ has positive Lebesgue measure on }\CC\}.
\end{equation*}
Then $\frac{\partial f}{\partial\bar z}=0$ on $S$. In particular, $f\in R(S)$.
\end{theorem}
Here, $\DD(z_0,r)=\{z\in\CC: |z-z_0|<r\}$ and $R(S)$ is the set of all continuous functions $g:K\to\CC$ which can be uniformly approximated (in sup-norm) on $S$ by functions holomorphic in a neighborhood of $S$. Edlund mentioned that he was unable to prove this for
 functions $f$ that are merely $C^1$. Also Coupet wrote about problems with $C^1$. Generally, the problems with $C^1$ stem from the fact
that the Bishop construction of a family of discs yields at best a $C^{k-\eps}$-regular family of boundaries of the disc, where $k$ is the regularity
of $f$. Therefore, if $f$ is $C^1$  we have no differentiable family of boundary curves in $M$, and we are unable to use arguments involving the  Jacobian of a mapping or the implicit function theorem for showing a diffeomorphism.

The main purpose of the present paper is to show that for $C^1$-functions, we still can prove a similar result for a smaller class of sets $S$. However the class is large enough to include the  so-called {\em fine domains}. See Section \ref{fine} below, where we give a reformulation of our result in the language of fine domains and finely holomorphic functions. So, as a corollary of Theorem~\ref{thm:main1} we get
\begin{theorem}\label{thm:3}
Let $K\subset\CC$ be a compact set and let $f:\CC\to\CC$ be a $C^1$-function such that 
the graph of $f$ over $K$ is pluripolar in $\CC^2$.
Put 
\begin{equation*}
S=\{z\in K:\text{ the set }\CC\setminus K\text{ is thin at }z\}.
\end{equation*}
Then $\frac{\partial f}{\partial\bar z}=0$ on $\overline{S}$. In particular, $f\in R(\overline S)$.
\end{theorem}
We recall the definition of a thin set below.

\medskip\noindent
{\bf Acknowledgement.} Part of the paper was done while the first author was visiting Indiana University Mathematics Department
as a Fulbright Scholar. He thanks the Department for its warm hospitality, especially Eric Bedford and Norm Levenberg.

\section{Thinness, fine topology and fine holomorphy}\label{fine}

We say that a set $F\subset\CC$ is {\em thin} at a point $z_0\in\CC$
if there exists a subharmonic function $u$ on $\CC$ such that $u\le -1$ on $F$ and $u(z_0)>-1$. 

Cartan already observed that when we provide $\CC$ with the {\em fine topology}, i.e. the coarsest topology that makes all subharmonic functions continuous, $F$ is thin at $a$ can be expressed as $a$ is not in the fine closure of $F$. It is a simple observation that if $\Omega$ is finely open and $z\in\Omega$ then there exists a compact set $K\subset \Omega$ which is a fine neighborhood of $z$.

On fine domains one can introduce finely holomorphic functions, which have many properties much similar to holomorphic functions, cf.~\cite{Fu81}. 

A function $f$ on a fine domain $\Omega$ is called {\em finely holomorphic} on $\Omega$, if for every $z\in \Omega$ there exists a compact fine neighborhood $z\in K\subset\Omega$, such that $f$ is a uniform limit on $K$ of rational functions with poles off $K$. Equivalently, there exists a $C^1$-function $f^*$ on $\CC$ such that $f=f^*$ on $K$ and $\bar\partial f^*=0$ on $K$.

Thus we can reformulate Theorem \ref{thm:3} in the language of fine holomorphy as follows.

\begin{theorem}\label{thm:4} Let $D\subset\CC$ be a finely open set and let $f:D\to\CC$ be a  function. Then the following conditions are equivalent:
\begin{enumerate}
\item $f$ is finely holomorphic on $D$;
\item $\Gamma_f(D)$ is pluripolar and for any point $z_0\in D$
there exist a compact fine neighborhood $K$ of $z_0$ and a $C^1$-function $f_K:\CC\to\CC$ such that
$f_K=f$ on $K$.
\end{enumerate}
\end{theorem}
\begin{proof} The implication $(1)\implies(2)$ was proved in the paper \cite{Edi-Marz-Wie}. The other implication
follows directly from Theorem~\ref{thm:3} and the definition of finely holomorphic function.
\end{proof}

\section{Proof of Theorems~\ref{thm:main1} and \ref{thm:3}}
Recall that a submanifold $M$ of $\CC^n$ is maximal totally real if  it has real dimension $n$ and
the tangent space $T_zM$ at any point $z\in M$ does no contain a complex line.

\begin{proof}[Proof of Theorem~\ref{thm:main1}]
Fix $z_0\in M$. Without loss of generality we may assume that $z_0=0\in\CC^n$.

According to E. Bishop (see \cite{Bishop}) because $M$ is maximal totally real we may write
in some neighborhood of the origin
$M=\{(x_1+ih_1(x_1,\dots,x_n),\dots, x_n+ih_n(x_1,\dots,x_n)\}$, where
$h_1,\dots,h_n$ are $C^1$ functions in a neighborhood of $0$ such that
$h_1(0)=\dots=h_n(0)=0$ and $\frac{\partial h_j}{\partial x_k}(0,0)=0$
for $j,k=1,\dots,n$. 

We recall the principal result from Coupet's paper \cite{Coupet}. Let $\DD$ be the unit disc in $\CC$ and let $\TT$ be its boundary. We denote by $T:L^p(\TT)\to L^p(\TT)$, $p\ge 1$, the Hilbert transform (or, harmonic conjugate transform).
In addition, for an $L^1$ function $\psi$ on $\TT$ we denote its harmonic extension by
\begin{equation}
\tilde\psi(z)=\frac1{2\pi}\int_{0}^{2\pi}\psi(e^{i\theta})P(e^{i\theta},z)d\theta,
\end{equation}
where $P(e^{i\theta},z)=\frac{1-|z|^2}{|e^{i\theta}-z|^2}$ is the Poisson kernel.
In particular, for $h\in L^1(\TT)$, 
\begin{equation}\label{jan2}
\tilde\psi +i\widetilde{T(\psi)}
\end{equation}
is a holomorphic function on $\DD$.

Recall the following result of B.~Coupet (see Th\'eor\`eme 1 in \cite{Coupet}).
\begin{thm}\label{theoreme:1}
Assume that $p>2n+1$. Then there exists a constant $\delta_0>0$ not depending on $h,n$ and $p$
such that for any function $k\in C^1(\TT\times\RR^{2n};\RR^n)$ with compact support and with
$\|k\|_{W^{1,p}(\TT\times\RR^{2n})}\le\delta_0$, the equation $u=T(h\circ u)+k$ has a unique
solution $u\in W^{1,p}(\TT\times\RR^{2n})$. 

Moreover, the harmonic extensions of $u$ and $h\circ u$ are of class $C^1$ on $\DD\times\RR^{2n}$ and
the mapping $\RR^{2n}\ni w\mapsto u(\cdot,w)\in W^{1,p}(\TT)$ is continuous.
\end{thm}

Fix $\phi\in C^\infty(\TT)$ such that $\phi=0$ on
$\TT^+=\{e^{i\theta}: |\theta|\le\frac\pi2\}$ and $\phi<0$ on
$\TT\setminus\TT^+$.
By Theorem~\ref{theoreme:1} we get (see also Part II in \cite{Coupet}) that there exist a small ball $B$ in $\RR^n$ centered at the origin and a unique function $u=u(z,\zeta,\xi):\ \overline \DD\times B\times B\to \RR^n$ such that
\begin{enumerate}
\item $u$ is harmonic in the first variable;
\item $u\in C(\bar\DD\times B\times B)\cap W^{1,p}(\bar\DD\times B\times B)$ for any fixed $p>2n+1$;
\item $u\in C^1(\DD\times B\times B)$;
\item $u(e^{i\theta},\zeta,\xi)=-T(h(u(\cdot,\zeta,\xi)))(e^{i\theta})+\zeta-\xi T(\phi)(e^{i\theta})$.
\end{enumerate}
Note that $u(e^{i\theta},\zeta,0)=\zeta$ and $u(0,\zeta,\xi)=\zeta$.
Define $H:\ \TT\times B\times B\to\CC^n$ by
\begin{equation}\label{jan1}
H(e^{i\theta},\zeta,\xi)=u(e^{i\theta},\zeta,\xi)+i\Big(h(u(e^{i\theta},\zeta,\xi))+\xi\phi(e^{i\theta})\Big)=\left(u+iTu\right)(e^{i\theta},\zeta,\xi),
\end{equation}
where (3) was used in the last equality.
Then for $\zeta,\xi\in B$, $H(e^{i\theta},\zeta,\xi)\in M$ for $|\theta|<\pi/2$. Moreover,
for fixed $\zeta,\xi$ the function $\tilde H (\cdot,\zeta,\xi)$ is holomorphic. 

Define $H:\ \TT\times B\times B\to\CC^n$ by
\begin{equation}\label{jan11}
H(e^{i\theta},\zeta,\xi)=u(e^{i\theta},\zeta,\xi)+i\Big(h(u(e^{i\theta},\zeta,\xi))+\xi\phi(e^{i\theta})\Big)=\left(u+iTu\right)(e^{i\theta},\zeta,\xi),
\end{equation}
where (3) was used in the last equality.
Then for $\zeta,\xi\in B$, $H(e^{i\theta},\zeta,\xi)\in M$ for $|\theta|<\pi/2$. Moreover,
for fixed $\zeta,\xi$ the function $\tilde H (\cdot,\zeta,\xi)$ is holomorphic. 
We have
\begin{equation}
\tilde H(0,\zeta,\xi)=\zeta+
\frac{i}{2\pi}\int_{0}^{2\pi}h(u(e^{i\theta},\zeta,\xi))d\theta+
i\xi\tilde\phi(0).
\end{equation}
Because $h(x)=o(x)$, we have $\frac{\partial \tilde H_j}{\partial\zeta_k}(0,0,0)=\delta_{jk}$ and 
$\frac{\partial \tilde H_j}{\partial\xi_k}(0,0,0)=\delta_{jk}i\tilde\phi(0)$, where $\delta_{jk}=1$ for $j=k$ and
$\delta_{jk}=0$ for $j\ne k$.


Now we will assume that $M\setminus P$ is thin at $0$ and show that this leads to a contradiction. Put $L=M\setminus (P\cup\{0\})$. Fix a neighborhood $W\subset\CC^n$ of $0$ and a negative plurisubharmonic function 
$v$ in $W$ such that $v\le-1$ on $W\cap L$ and $v(0)>-\frac{1}{4}$. 
We put $ S'=\{z\in\CC^n: v(z)>-\frac13\}$. Note that $0\in S'$.

Fix also a plurisubharmonic function $U$ on $\CC^n$ such that $U=-\infty$ on $P$.
Then for fixed $\zeta,\xi$ the function $u(w,\zeta,\xi)=U(H(w,\zeta,\xi))$ is subharmonic on $\DD$.

We are interested in the set $\{(\theta,\zeta,\xi): H(e^{i\theta},\zeta,\xi)\not\in L\}$.

\begin{lemma} Let $\zeta,\xi$ be near zero such that $ H(\overline\DD,\zeta,\xi)\subset W$
and that $H(0,\zeta,\xi)\in S'$. Then $\{\theta: |\theta|<\frac{\pi}{2}, H(e^{i\theta},\zeta,\xi)\not\in L\}$
is of positive measure.
\end{lemma}

\begin{proof} We have
\begin{equation}
v(H(0,\zeta,\xi))\le\frac{1}{2\pi}\int_{-\frac\pi2}^{\frac\pi2} v( H(e^{i\theta},\zeta,\xi))d\theta.
\end{equation}
If for almost all $\theta\in(-\frac\pi2,\frac\pi2)$ we have $H(e^{i\theta},\zeta,\xi)\in L$ 
then
\begin{equation}
-\frac13<v(H(0,\zeta,\xi))\le-\frac{1}{2}.
\end{equation}
A contradiction.
\end{proof}

\begin{lemma} For a set of positive measure $A\subset\RR^{2n}$ we have
$H(0,\zeta,\xi)\in S'$, $(\zeta,\xi)\in A$.
\end{lemma}
\begin{proof}
Indeed, $\det[\frac{\partial H}{\partial\zeta}](0,0,0)=i$. Hence, for small
$\zeta, \xi\in\RR^n$ we have
$\det[\frac{\partial H}{\partial\zeta}](0,\zeta,\xi)\not=0$.
Note that for any fixed $\xi$ near zero the set $(H)^{-1}(S')$ is of positive
measure.
\end{proof}

The set $\{H(0,\zeta,\xi): (\zeta,\xi)\in A\}$ is of positive measure in $\CC^n$, because $H$ is $C^1$ at the origin.
Then we get that $U\equiv-\infty$ on $\CC^n$. A contradiction.
\end{proof}

\begin{proof}[Proof of Theorem~\ref{thm:3}]
Fix a point $z_0\in S$. We want to show that $\frac{\partial f}{\partial\overline z}(z_0)=0$.
Assume that this is not the case. Then in some neighborhood of $z_0$ we have
\begin{equation}\label{eq:12}
\frac{\partial f}{\partial\overline z}\not=0. 
\end{equation}
We put $M=\{(z,f(z)):z\text{ near } z_0\}$. From \eqref{eq:12} we see that $M$ is a maximal totally real submanifold.
This would contradict Theorem~\ref{thm:main1}, where $P=\{(z,f(z)): z\in K\}$ is the pluripolar set.
\end{proof}

\section{An example}\label{example}

The following example indicates how the compact sets that we considered may come about.

\begin{example}
Take sequences $a_n\searrow0$, $r_n>0$, and take $K=\{0\}\cup\cup_{n=1}^\infty\overline{\DD(a_n,r_n)}$.
Note that for any $C^1$-function $f:\CC\to\CC$ if $\Gamma_f(K)$ is pluripolar then by Shcherbina's result
$\frac{\partial f}{\partial\bar z}(0)=0$. Now, take $r_n$ so small that all discs $\DD(a_n,r_n)$ are disjoint.

Fix one of these discs, say $\DD(a_m,r_m)$. We take a dense subset $\{b_n\}$ in this disc. If we carefully remove balls
about points $b_n$ then we may obtain a set $L_m$, so that $\CC\setminus L_m$ is thin at $a_m$.
We do this for all discs and we get a new compact set $\tilde K=\{0\}\cup\cup_{m=1}^\infty L_m$ such that
its interior is empty. Again take a $C^1$ function $f:\CC\to\CC$ such that $\Gamma_f(\tilde K)$ is pluripolar.
Then neither Shcherbina's nor Edlund's result can be applied. But our main theorem gives $\frac{\partial f}{\partial\bar z}(0)=0$. The reader can find similar, but even more sophisticated examples of this type.
\end{example}

%
\bibliographystyle{amsplain}


\end{document}